\documentclass[a4paper]{amsart}

\usepackage{amsmath}
\usepackage{amssymb}
\usepackage{color}
\usepackage{graphicx}
\usepackage{hyperref}

\newtheorem{theorem}{Theorem}
\newtheorem*{theorem*}{Theorem}
\newtheorem{lemma}{Lemma}
\newtheorem{lemma_and_definition}[lemma]{Lemma (and Definition)}

\newtheorem*{definition*}{Definition}

\newcommand{\R}{\mathbb{R}}
\newcommand{\RP}{\R\mathrm{P}}
\newcommand{\Orth}{\mathit{O}}
\newcommand{\GL}{\mathit{GL}\,}

\title{A unique representation of polyhedral types.
  Centering via M\"obius transformations}

\author{Boris A.~Springborn}

\thanks{The author is supported by DFG Research Center
    \textsc{Matheon}.}

\date{16 March 2004}

\begin{document}

\begin{abstract}
  For $n\geq 3$ distinct points in the $d$-dimensional unit sphere
  $S^d\subset\R^{d+1}$, there exists a M\"obius transformation such that the
  barycenter of the transformed points is the origin. This M\"obius
  transformation is unique up to post-composition by a rotation. We prove
  this lemma and apply it to prove the uniqueness part of a representation
  theorem for 3-dimensional polytopes as claimed by Ziegler (1995): For each
  polyhedral type there is a unique representative (up to isometry) with
  edges tangent to the unit sphere such that the origin is the barycenter of
  the points where the edges touch the sphere.  
\end{abstract}

\maketitle 


In today's language, Steinitz' {\em fundamental theorem of convex
  types}\/~\cite{steinitz22:_polyed_raumein},
\cite{steinitz34:_vorles_theor_polyed} (for a modern treatment
see~\cite{gruenbaum03:_convex_polyt}, \cite{ziegler95:_lectur_polyt}) states
that the combinatorial types of convex 3-dimensional polyhedra correspond to
the strongly regular cell decompositions of the 2-sphere. (A cell complex is
{\em regular} if the closed cells are attached without identifications on the
boundary. A regular cell complex is {\em strongly regular}\/ if the
intersection of two closed cells is a closed cell or empty.) 

Gr\"unbaum and Shephard~\cite{gruenbaum87:_some_probl_on_polyh} posed the
question whether for every combinatorial type there is a polyhedron with
edges tangent to a sphere.  This question has been answered affirmatively:

\begin{theorem}[Koebe~\cite{koebe36:_kontak_abbil},
  Andreev~\cite{andreev70a}\cite{andreev70b},
  Thurston~\cite{thurston:_geomet_topol_three_manif}, Brightwell and
  Scheinerman~\cite{brightwell93:_repres_of_planar_graph},
  Schramm~\cite{schramm92:_how}]
  \label{thm:koebe}
  For every combinatorial type of convex $3$-dimensional polyhedra, there is
  a representative with edges tangent to the unit sphere. This representative
  is unique up to projective transformations which fix the sphere and do not
  make the polyhedron intersect the plane at infinity.
\end{theorem}

A proof which makes use of a variational principle was given by
A.~Bo\-ben\-ko and the author~\cite{bobenko04:_variat},
\cite{springborn03:_variat}.

The purpose of this article is to prove Theorem~\ref{thm:main} below,
which singles out a unique representative for each convex type. (The proof
given here is also contained in the author's doctoral
dissertation~\cite{springborn03:_variat}.)  The claim of
Theorem~\ref{thm:main} is not new (see
Ziegler~\cite{ziegler95:_lectur_polyt}, p.~118, and the second edition of
Gr\"unbaum's classic~\cite{gruenbaum03:_convex_polyt}, p.~296a) but this
proof seems to be.

\begin{theorem}\label{thm:main}
  For every combinatorial type of convex 3-dimensional polyhedra there is a
  unique polyhedron (up to isometry) with edges tangent to the unit sphere
  $S^2\subset\R^3$, such that the origin\/ $0\in\R^3$ is the barycenter of
  the points where the edges touch the sphere.
\end{theorem}

Theorem~\ref{thm:main} follows from Theorem~\ref{thm:koebe} and
Lemma~\ref{lem:moebius_center} below (with $d=2$). Indeed, the projective
transformations of $\RP^{d+1}$ that fix $S^d$ correspond to the M\"obius
transformations of $S^d$. Lemma~\ref{lem:moebius_center} is also of interest
in its own right.

\begin{lemma}
  \label{lem:moebius_center}
  Let $v_1,\ldots,v_n$ be $n\geq 3$ distinct points in the $d$-dimensional
  unit sphere $S^d\subset\R^{d+1}$. There exits a M\"obius transformation $T$
  of $S^d$, such that
  \begin{equation*}
    \sum_{j=1}^n Tv_j=0.
  \end{equation*}
  If $\widetilde T$ is another such M\"obius transformation, then
  $\widetilde{T}=RT$, where $R$ is an isometry of $S^d$.
\end{lemma}

Our proof of Lemma~\ref{lem:moebius_center} is based on the fundamental
relationship between projective, hyperbolic, and M\"obius geometry. The
equation
 \begin{equation*}
   -x_0^2+x_1^2+x_2^2+\ldots+x_{d+1}^2=0
 \end{equation*}
represents the $d$-dimensional sphere $S^d$ as a quadric in
$(d+1)$-dimensional projective space $\RP^{d+1}$. The group of projective
transformations of $\RP^{d+1}$ which fix $S^d$ is $\Orth(d+1,1)/\{\pm 1\}$,
where the orthogonal group $\Orth(d+1,1)\subset\GL(d+2)$ acts linearly on the
homogeneous coordinates.  At the same time, $\Orth(d+1,1)/\{\pm 1\}$ acts
faithfully as the M\"obius group on $S^d$, and as the isometry group of
$(d+1)$-dimensional hyperbolic space $H^{d+1}$, which is identified with the
open ball bounded by $S^d$ (the Klein model of hyperbolic space). For a
detailed account of this classical material see, for example,
Hertrich-Jeromin~\cite{hertrich-jeromin03:_introd_moebius_differ_geomet} and
Kulkarni, Pinkall~\cite{kulkarni88:_confor_geomet}.

A similar interplay of geometries leads Bern and Eppstein, to another choice
of a unique representative for each polyhedral type. Given $n$ spheres in
$S^d$, Bern and Eppstein apply that M\"obius transformation which makes the
smallest sphere as large as possible~\cite{bern01:_optim_moebius}. It is not
difficult to see that this M\"obius transformation is unique up to
post-composition with a rotation if $n\geq 3$. Since edge-tangent polyhedra
correspond to circle packings, this leads to another choice of unique
representative for each polyhedral
type~\cite{eppstein03:_hyper_geomet_moebius_trans_geomet_optim}.

For symmetric polyhedral types (more precisely, for those polyhedral types
with a symmetry group of orientation preserving isomorphisms which is not
just a cyclic group) the unique representative of Bern and Eppstein coincides
with ours.

\section*{Proof of Lemma~\ref{lem:moebius_center}.}

The M\"obius transformations of $S^d$ correspond to isometries of the
hyperbolic space $H^{d+1}$, of which $S^d$ is the infinite boundary.

For $n\geq 3$ points $v_1,\ldots,v_n\in S^d$, we are going show
that there is a unique point $x\in H^{d+1}$ such that the sum of the
`distances' to $v_1,\ldots,v_n$ is minimal. Of course, the distance to a
point in the infinite boundary is infinite. The right quantity to use instead
is the distance to a horosphere through the infinite point (see the figure).

\begin{figure}%
\hfill%
\input{horosphere_ball.pstex_t}%
\hfill%
\input{horosphere_halfspace.pstex_t}%
\hspace*{\fill} \\[\baselineskip]%
\parbox{0.9\textwidth}{%
{\sc Figure.} The `distance' to an infinite point $v$ is measured by cutting 
off at some horosphere through $v$. {\em Left:}\/ Poincar\'e ball model.
{\em Right:}\/ half-space model.}
\end{figure}

\begin{definition*}
  For a horosphere $h$\/ in $H^{d+1}$, define
  \begin{gather*}
    \delta_h:H^{d+1}\rightarrow \R,\\
    \delta_h(x) = 
    \begin{cases}
      -\operatorname{dist}(x, h) & \text{if $x$ is inside $h$}, \\
      0 & \text{if $x\in h$}, \\
      \operatorname{dist}(x, h) & \text{if $x$ is outside $h$},
    \end{cases}
  \end{gather*}
  where $\operatorname{dist}(x, h)$ is the distance from the point $x$ to the
  horosphere $h$.
\end{definition*}

Suppose $v$ is the infinite point of the horosphere $h$. Then the shortest
path from $x$ to $h$ lies on the geodesic connecting $x$ and $v$. If $h'$ is
another horosphere through $v$, then $\delta_h-\delta_{h'}$ is constant. If
$g:\R\rightarrow H^{d+1}$ is an arc-length parametrized geodesic, then
$\delta_h\circ g$ is a strictly convex function, unless $v$ is an infinite
endpoint of the geodesic $g$. In that case, $\delta_h\circ g(s)=\pm(s-s_0)$.
These claims are straightforward to prove using the Poincar\'e half-space
model, where hyperbolic space is identified with the upper half space:
\begin{equation*}
  H^{d+1}=\big\{(x_0,\ldots,x_d)\in\R^{d+1}\;\big|\;x_0>0\big\},
\end{equation*}
and the metric is
\begin{equation*}
  ds^2=\frac{1}{x_0^2}\big(dx_0^2+dx_1^2+\cdots+dx_d^2\big)
\end{equation*}
Also, one finds that, as $x\in H^{d+1}$ approaches the infinite boundary,
\begin{equation*}
  \lim_{x\rightarrow\infty}\sum_{j=1}^n \delta_{h_j}(x)=\infty,
\end{equation*}
where $h_j$ are horospheres through different infinite points and $n\geq3$.
Thus, the following definition of the {\em point of minimal distance sum} is
proper.

\begin{lemma_and_definition}\hspace{0.5em}
  Let $v_1,\ldots,v_n$ be $n$ points in the infinite boundary of $H^{d+1}$,
  where $n\geq 3$. Choose horospheres $h_1,\ldots,h_n$ through
  $v_1,\ldots,v_n$, respectively. Then there is a unique point $x\in H^{d+1}$
  for which $\sum_{j=1}^n \delta_{h_j}(x)$ is minimal. This point $x$ does not
  depend on the choice of horospheres. It is the {\em{}point of minimal
    distance sum}\/ from the infinite points $v_1,\ldots,v_n$.
\end{lemma_and_definition}

In the Poincar\'e ball model, hyperbolic space is identified with the unit
ball as in the Klein model, but the metric is
$  
ds^2=\frac{4}{(1-\sum x_j^2)_{\rule{0pt}{1ex}}^2}\sum \,dx_j^2.
$
(Since the Klein model and the Poincar\'e ball model agree on the infinite
boundary and in the center of the sphere, one might as well use the Klein
model in the following lemma.)

\begin{lemma}
  Let $v_1,\ldots,v_n$ be $n\geq3$ different points in the infinite boundary
  of $H^{d+1}$. In the Poincar\'e ball model, $v_j\in S^d\subset\R^{d+1}$.
  The origin is the point of minimal distance sum, if and only if $\sum
  v_j=0$.
\end{lemma}

\begin{proof}
  If $h_j$ is a horosphere through $v_j$, then the gradient of $\delta_{h_j}$
  at the origin is the unit vector $-\frac{1}{2}v_j$. 
\end{proof}

Lemma~\ref{lem:moebius_center} is now almost
immediate. Let $x$ be the point of minimal distance sum from the
$v_1,\ldots,v_n$ in the Poincar\'e ball model.  There is a hyperbolic
isometry $T$ which moves $x$ into the origin. If $\widetilde T$ is another
hyperbolic isometry which moves $x$ into the origin, then $\widetilde T=RT$,
with $R$ is an orthogonal transformation of $\R^{d+1}$.

This concludes the proof of Lemma~\ref{lem:moebius_center}.

\section*{Acknowledgements}

I would like to thank Alexander Bobenko and G\"unter Ziegler for making me
familiar with the problem of finding unique representatives for polyhedral
types, and Ulrich Pinkall, who has provided the essential insight for this
solution. 

\bibliographystyle{plain}
\bibliography{unirep}

\end{document}